\documentclass[11pt]{amsart}
\usepackage{amsmath,amssymb,standardmath}
\usepackage[margin=1in]{geometry}

\newcommand{\hook}[3]{\chi^{(#1|#2)}_{#3}}
\newcommand{\hookGS}[1]{F_{#1}}
\newcommand{\RGS}[1]{R_{#1}}
\newcommand{\cx}[2]{c^{#1}_{#2}}
\renewcommand{\r}[1]{r^{#1}}

\title{On the Number of Factorizations of a Full Cycle}
\author{John Irving}
\email{jcirving@math.uwaterloo.ca}
\address{Department of Combinatorics \& Optimization, University of Waterloo, Canada}
\date{\today}

\begin{document}
\begin{abstract}
We give a new expression for the number of factorizations of a full cycle into an ordered product of permutations
of specified cycle types. This is done through purely algebraic means, extending work of Biane.  We deduce from
our result a formula of Poulalhon and Schaeffer that was previously derived through an intricate combinatorial
argument.
\end{abstract}

\maketitle

\section{Notation}

Our notation is generally consistent with Macdonald~\cite{bk:macdonald}. We write $\l \ptn n$ (or $|\l|=n$) and
$\len{\l}=k$ to indicate that $\l$ is a partition of $n$ into $k$ parts; that is, $\l=(\l_1,\ldots,\l_k)$ with
$\l_1 \geq \cdots \geq \l_k \geq 1$ and $\l_1+\ldots+\l_k=n$. If $\l$ has exactly $m_i$ parts equal to $i$ then we
write $\l = [1^{m_1} 2^{m_2} \cdots]$, suppressing terms with $m_i=0$. We also define $\z{\l} = \prod_i i^{m_i}
m_i!$ and $\aut{\l} = \prod_i m_i!$. A \emph{hook} is a partition of the form $[1^b,\,a+1]$ with $a, b \geq 0$. We
use Frobenius notation for hooks, writing $(a|b)$ in place of $[1^b,\,a+1]$.

The conjugacy class of the symmetric group $\Sym{n}$ consisting of all $n!/\z{\l}$ permutations of cycle type $\l
\ptn n$ will be denoted by $\Class{\l}$. The irreducible characters $\chi^{\l}$ of $\Sym{n}$ are naturally indexed
by partitions $\l$ of $n$, and we use the usual notation $\ch{\l}{\mu}$ for the common value of $\chi^\l$ at any
element of $\Class{\mu}$. We write $\f[\l]$ for the degree $\ch{\l}{[1^n]}$ of $\chi^\l$.

For vectors $\mathbf{j}=(j_1,\ldots,j_m)$ and $\mathbf{x}=(x_1,\ldots,x_m)$ we use the abbreviations $\mathbf{j!}
= j_1!\cdots j_m!$ and $\mathbf{x}^\mathbf{j} = x_1^{j_1} \cdots x_m^{j_m}$. Finally, if $\a \in \Q$ and $f \in
\Q[[\mathbf{x}]]$ is a formal power series, then we write $[\a \mathbf{x}^{\mathbf{j}}]\,f(\mathbf{x})$ for the
coefficient of the monomial $\a \mathbf{x}^{\mathbf{j}}$ in $f(\mathbf{x})$.

\section{Factorizations of Full Cycles}

Given $\l, \a_1,\ldots,\a_m \ptn n$, let $\cx{\l}{\a_1,\ldots,\a_m}$ be the number of factorizations in $\Sym{n}$
of a given permutation $\p \in \Class{\l}$ as an ordered product $\p = \s_1 \ldots \s_m$, with $\s_i \in
\Class{\a_i}$ for all $i$. The problem of evaluating $\cx{\l}{\a_1,\ldots,\a_m}$ for various $\l$ and $\a_i$ has
attracted a good deal of attention and is linked to various questions in algebra,  geometry, and physics. For
details on the history of this problem and its connections to other areas of mathematics, we direct the reader
to~\cite{goupil-schaeffer:largecycles} and the references therein. Here we focus on the particularly well-studied
case $\l = (n)$, which corresponds to counting factorizations of the full cycle $(1\,2\,\cdots\,n) \in \Sym{n}$
into factors of specified cycle types.

While it is straightforward to express $\smash{\cx{(n)}{\a_1,\ldots,\a_m}}$ as a character sum
(see~\eqref{eq:gencx2} below), the appearance of alternating signs in this sum --- and the resulting cancellations
--- preclude asymptotic analysis. Goupil and Schaeffer~\cite{goupil-schaeffer:largecycles} overcame this difficulty
in the case $m=2$ by interpreting certain characters combinatorially (viz. the Murnaghan-Nakayama rule) and
employing a sequence of bijections in which a sign-reversing involution accounts for cancellations. This leads to
an expression for $\smash{\cx{(n)}{\a,\b}}$ as a sum of \emph{positive} terms, which in turn permits nontrivial
asymptotics. Poulalhon and Schaeffer~\cite{poulalhon-schaeffer:largecycles} later extended this argument to arrive
at a similar formula for $\smash{\cx{(n)}{\a_1,\ldots,\a_m}}$.

Biane~\cite{biane:largecycles} has recently given a remarkably succinct algebraic derivation of Goupil and
Schaeffer's formula for $\smash{\cx{(n)}{\a,\b}}$.  Our purpose here is to extend his method to give a new
expression for $\smash{\cx{(n)}{\a_1,\ldots,\a_m}}$ as a sum of positive contributions.  In particular, if for $\g
=(\g_1,\g_2,\ldots) \ptn n$ we define the polynomial $\RGS{\g}(x,y)$ and the \emph{nonnegative} constants
$\r{\g}_{j,k}$ by
\begin{equation} \label{eq:RGS}
    \RGS{\g}(x,y)
        := \frac{1}{2y} \prod_{i \geq 1} ((x+y)^{\g_i}-(x-y)^{\g_i})
        =  \sum_{j+k=n-1} \r{\g}_{j,k} x^{j} y^k,
\end{equation}
then our main result is the following:

\begin{thm}
\label{thm:main} Let $\a_1,\ldots,\a_m \ptn n$ and, for $\l=[1^{m_1} 2^{m_2} 3^{m_3} \cdots]$, let $2\l-1=[1^{m_1}
3^{m_2} 5^{m_3} \cdots]$. Set $\mathbf{x}=(x_1,\ldots,x_m)$ and let $e_{\l}(\mathbf{x})$ denote the elementary
symmetric function in $x_1,\ldots,x_m$ indexed by $\l$.  Then
\begin{equation*}
        \cx{(n)}{\a_1,\ldots,\a_m}
    = \frac{n^{m-1}}{2^{(n-1)(m-1)} \prod_i \z{\a_i}}  \sum_{\mathbf{j+k=n-1}}
        [\mathbf{x}^{\mathbf{j}}] \!\!\sum_{\len{\l}=n-1} \!\!\!\frac{e_{2\l-1}(\mathbf{x})}{\aut{\l}}
         \cdot \prod_{i=1}^m j_i!\,k_i!\,\r{\a_i}_{j_i,k_i},
\end{equation*}
where the outer sum extends over all vectors $\mathbf{j}=(j_1,\ldots,j_m)$ and $\mathbf{k}=(k_1,\ldots,k_m)$ of
nonnegative integers such that $j_i+k_i=n-1$ for all $i$, and the inner sum over all partitions $\l$ with $n-1$
parts.
\end{thm}

A proof of Theorem~\ref{thm:main} is given in the next section. In \S4, we use this result to deduce Poulalhon and
Schaeffer's formula for $\smash{\cx{(n)}{\a_1,\ldots,\a_m}}$ (listed here as Theorem~\ref{thm:ps}), thereby giving
a purely algebraic derivation that avoids the detailed combinatorial constructions
in~\cite{poulalhon-schaeffer:largecycles}.

\section{Proof of the Main Result}

\newcommand{\dmu}[1]{d\mu(#1)}
\newcommand{\tprod}{{\textstyle \prod}}

It is well known that the class sums $\mathsf{K}_{\l} = \sum_{\s \in \Class{\l}} \s$ (for $\l \ptn n$) form a
basis of the centre of the group algebra $\C\Sym{n}$. Indeed, the linearization relations
$\mathsf{K}_{\a_1}\cdots\mathsf{K}_{\a_m} = \sum_{\l \ptn n} \cx{\l}{\a_1,\ldots,\a_m} \mathsf{K}_{\l}$ identify
the constants $\cx{\l}{\a_1,\ldots,\a_m}$ as the connection coefficients of $\C\Sym{n}$. Using character theory
(in particular, by expressing $\mathsf{K}_{\l}$ in terms of the central primitive idempotents of $\C\Sym{n}$) one
finds that
\begin{equation*}\label{eq:gencx}
    \cx{\l}{\a_1,\ldots,\a_m}
    = \frac{n!^{m-1}}{\z{\a_1} \cdots \z{\a_m}} \sum_{\b \ptn n} \frac{\ch{\b}{\a_1} \cdots \ch{\b}{\a_m}}
                                {(\f[\b])^{m-1}} \ch{\b}{\l}.
\end{equation*}
This sum is generally intractable but simplifies considerably in the case $\l = (n)$, since here $\ch{\b}{\l}$
vanishes when $\b$ is not a hook; in particular, the Murnaghan-Nakayama rule~\cite{bk:stanley2} implies
$\ch{\b}{(n)} = (-1)^b$ if $\b = (a | b)$, while $\ch{\b}{(n)}=0$ otherwise. Moreover, the hook-length formula
gives $\f[(a|b)] = \binom{a+b}{b}$, so
\begin{align}\label{eq:gencx2}
    \cx{(n)}{\a_1,\ldots,\a_m}
    = \frac{n^{m-1}}{\z{\a_1} \cdots \z{\a_m}}
     \sum_{a+b=n-1} (a!\,b!)^{m-1} \ch{(a|b)}{\a_1} \cdots \ch{(a|b)}{\a_m} (-1)^b.
\end{align}

Let $\mu$ be the measure on $\mathbb{C}$ defined by the density $ \dmu{z} = \tfrac{1}{\pi} e^{-|z|^2} dz, $ where
$dz$ is the standard Lebesgue density (\emph{i.e.} $dz = ds\,dt$ for $z=s + t\sqrt{-1}$). Following
Biane~\cite{biane:largecycles}, we shall make use of the formula
\begin{equation}
\label{eq:wick}
    \int_{\mathbb{C}} z^j \bar{z}^k \dmu{z} = j! \,\delta_{jk},
\end{equation}
which is easily verified by changing to polar form.

\noindent \textbf{Proof of Theorem~\ref{thm:main}}:
 For $\g  \ptn n$, let $\hookGS{\g}(u,v) = \sum \hook{a}{b}{\g}
u^a v^b$ be the generating series for hook characters, where the sum extends over all pairs $(a,b)$ of nonnegative
integers with $a+b=n-1$. Then
\begin{align*}
    \frac{1}{(n-1)!}&(u_1\cdots u_m - v_1 \cdots v_m)^{n-1} \prod_{i=1}^m
    F_{\a_i}(\bar{u}_i,\bar{v}_i) \notag \\
    &= \sum_{a+b=n-1} \frac{u_1^a \cdots u_m^a \cdot v_1^b \cdots v_m^b}{a!\,b!} (-1)^b
         \,\prod_{i=1}^m  \sum_{a_i + b_i = n-1} \hook{a_i}{b_i}{\a_i} \bar{u}_i^{a_i} \bar{v}_i^{b_i}.
         \label{eq:bigseries}
\end{align*}
Consider the effect of integrating the RHS with respect to $\dmu{\mathbf{u},\mathbf{v}}:=\prod_{i=1}^m \dmu{u_i}
\dmu{v_i}$. Using~\eqref{eq:wick}, note that all monomials $\frac{(-1)^b}{a!b!} \prod_i \hook{a_i}{b_i}{a_i} u_i^a
\bar{u}_i^{a_i} v_i^b \bar{v}_i^{b_i}$  vanish except those with $a_i=a$ and $b_i=b$ for all $i$, and each
monomial of this special form is replaced by $(-1)^b (a!\,b!)^{m-1} \prod_i \hook{a}{b}{\a_i}$.  Thus we obtain
\begin{align*}
    \int_{\C^{2m}} (u_1\cdots u_m - v_1 \cdots &v_m)^{n-1} \prod_{i=1}^m
        F_{\a_i}(\bar{u}_i,\bar{v}_i) \,\dmu{\mathbf{u},\mathbf{v}} \\
        &=  (n-1)! \!\! \sum_{a+b=n-1} (a!\,b!)^{m-1} \ch{(a|b)}{\a_1} \cdots \ch{(a|b)}{\a_m} (-1)^b.
\end{align*}
Let $I$ be the integral on the LHS, and change variables by letting $ u_i=(y_i+x_i)/\sqrt{2}$,
$v_i=(y_i-x_i)/\sqrt{2}$. As an immediate consequence of the Murnaghan-Nakayama rule we have
\begin{equation*}
\label{eq:hookGS}
    \hookGS{\g}(u,v)
    = \frac{1}{u+v} \prod_{i \geq 1} (u^{\g_i} - (-v)^{\g_i})
\end{equation*}
for a partition $\g = (\g_1,\g_2,\ldots)$. Thus~\eqref{eq:RGS} gives $\hookGS{\g}(y+x,y-x) = \RGS{\g}(x,y)$, and
since $\hookGS{\a_i}$ is homogeneous of degree $n-1$ the change of variables yields
$\hookGS{\a_i}(\bar{u}_i,\bar{v}_i) = 2^{-(n-1)/2}\RGS{\a_i}(\bar{x}_i,\bar{y}_i)$ for all $i$. Furthermore, it is
easy to check that $\dmu{\mathbf{u},\mathbf{v}} = \dmu{\mathbf{x},\mathbf{y}}$ and
\begin{align*}
    u_1\cdots u_m - v_1 \cdots v_m
    = \frac{1}{\sqrt{2^m}} \bigg(\prod_{i=1}^m (y_i+x_i) - \prod_{i=1}^m (y_i - x_i)\bigg)
    = \frac{2 y_1\cdots y_m}{\sqrt{2^m}} \sum_{s \geq 1} e_{2s-1}( \mathbf{x}/\mathbf{y}),
\end{align*}
where $\mathbf{x}/\mathbf{y}=(\tfrac{x_1}{y_1},\ldots,\tfrac{x_m}{y_m})$. Thus, with the aid of~\eqref{eq:wick},
we get
\begin{align*}
    I &= \frac{1}{2^{(n-1)(m-1)}} \int_{\C^{2m}} \bigg( y_1\cdots y_m \sum_{s \geq 1}
    e_{2s-1}(\mathbf{x}/\mathbf{y}) \bigg)^{n-1} \prod_{i=1}^m
        \RGS{\a_i}(\bar{x}_i,\bar{y}_i) \,\dmu{\mathbf{x},\mathbf{y}} \\
    &= \frac{1}{2^{(n-1)(m-1)}}  \sum_{\mathbf{j,k}} \, \mathbf{j}!\,\mathbf{k}!\,
    [\mathbf{x}^\mathbf{j} \mathbf{y}^\mathbf{k}]\, \bigg( y_1\cdots y_m \sum_{s \geq 1}
    e_{2s-1}(\mathbf{x}/\mathbf{y}) \bigg)^{n-1}\!\!
               \cdot \,[\bar{\mathbf{x}}^{\mathbf{j}} \bar{\mathbf{y}}^\mathbf{k}] \prod_{i=1}^m \RGS{\a_i}(\bar{x}_i,\bar{y}_i) \\
    &= \frac{1}{2^{(n-1)(m-1)}}  \sum_{\mathbf{j+k=n-1}} \, \mathbf{j}!\,\mathbf{k}!\,
    [\mathbf{x}^\mathbf{j}]\, \bigg( \sum_{s \geq 1}
    e_{2s-1}(\mathbf{x}) \bigg)^{n-1}
               \prod_{i=1}^m \r{\a_i}_{j_i,k_i} \\
    &= \frac{(n-1)!}{2^{(n-1)(m-1)}}  \sum_{\mathbf{j+k=n-1}} \, \mathbf{j}!\,\mathbf{k}!\,
    [\mathbf{x}^\mathbf{j}]\, \sum_{\len{\l} = n-1} \frac{e_{2\l-1}(\mathbf{x})}{\aut{\l}}
               \prod_{i=1}^m \r{\a_i}_{j_i,k_i}.
\end{align*}
The result now follows from~\eqref{eq:gencx2}. \qed

\section{Recovery of Poulalhon \& Schaeffer's Formula}

\newcommand{\SGS}[1]{S_{#1}}
\newcommand{\PGS}[2]{P_{#1}^{#2}}
\newcommand{\Dop}{\mathfrak{D}}

We require some extra notation to state the Poulalhon-Schaeffer formula for $\cx{(n)}{\a_1,\ldots,\a_m}$. First,
we define symmetric polynomials $\SGS{p}(x_1,\ldots,x_l)$ by setting $\SGS{0}(x_1,\ldots,x_l)=1$ and
$$
    \SGS{p}(x_1,\ldots,x_l) = \sum_{p_1+\cdots+p_l = p} \,\prod_{i=1}^l \frac{1}{x_i} \binom{x_i}{2p_i+1}
$$
for $p > 0$. Note that these have the simple generating series
\begin{equation} \label{eq:SGS}
    \sum_{p \geq 0} \SGS{p}(x_1,\ldots,x_l) t^{2p} = \prod_{i=1}^l \frac{(1+t)^{x_i}-(1-t)^{x_i}}{2x_i t},
\end{equation}
which is obviously closely related to our series $\RGS{\g}(x,y)$ (see~\eqref{eq:RGS}). We also introduce an
operator $\Dop$ on $\Q[[x_1,\ldots,x_m]]$ defined as follows: For each $i$ and all $j \geq 0$ set
$\Dop(x_i^j)=x_i(x_i-1)\cdots(x_i-j+1)$, and extend the action of $\Dop$ multiplicatively to monomials
$x_1^{j_1}\cdots x_m^{j_m}$ and then linearly to all of $\Q[[x_1,\ldots,x_m]]$.  Finally, we define polynomials
$\PGS{a}{b}(x_1,\ldots,x_m)$ by setting $\PGS{0}{b}(x_1,\ldots,x_m)=1$ for all $b \geq 1$ and letting
\begin{equation}\label{eq:PGS}
    \PGS{a}{b}(x_1,\ldots,x_m)
    \,= \!\!
    \sum_{\substack{\l \ptn a \\ \len{\l} \leq b}}
    \!\!\Dop\Big(\frac{e_{2\l+1}(\mathbf{x})}{\aut{\l}}\Big)
\end{equation}
for $a,b \geq 1$, where $2\l+1=[3^{m_1} 5^{m_2} \cdots]$ when $\l = [1^{m_1} 2^{m_2}\cdots]$. Then the main result
of~\cite{poulalhon-schaeffer:largecycles} is the  following formula for $\cx{(n)}{\a_1,\ldots,\a_m}$.

\begin{thm}\label{thm:ps}
Let $\a_1,\ldots,\a_m \ptn n$ and set $r_i=n-\len{\a_i}$ for all $i$.  Let $g=\frac{1}{2}(\sum_i r_i - n+1)$. If
$g$ is a nonnegative integer, then
\begin{equation*}
    \cx{(n)}{\a_1,\ldots,\a_m}
    = \frac{n^{m-1}}{2^{2g} \prod_{i} \aut{\a_i}}
        \sum
        \PGS{q}{n-1}(\mathbf{r-2p}) \prod_{i=1}^m (\len{\a_i}+2p_i-1)! \,\SGS{p_i}(\a_i),
\end{equation*}
where $\mathbf{r-2p} = (r_1-2p_1,\ldots,r_m-2p_m)$ and the sum extends over all tuples $(q,p_1,\ldots,p_m)$ of
nonnegative integers with $q+p_1+\cdots+p_m=g$.
\end{thm}

Before proceeding to deduce this result from Theorem~\ref{thm:main}, we pause for a few remarks. First, the
integer $g$ identified in Theorem~\ref{thm:ps} is called the \emph{genus} of the associated factorizations of
$(1\,2\,\cdots\,n)$, and it has well-understood geometric meaning; see~\cite{melou-schaeffer:constellations}, for
example.  The primary benefit of the Poulalhon-Schaeffer formula (over Theorem~\ref{thm:main}) is that the
dependence on genus is explicit.  For instance, when $g=0$  it is immediately clear that Theorem~\ref{thm:ps}
reduces to the very simple
$$
    \cx{(n)}{\a_1,\ldots,\a_m} = n^{m-1} \prod_{i=1}^m \frac{(\len{\a_i}-1)!}{\aut{\a_i}}.
$$
The $\smash{\cx{(n)}{\a_1,\ldots,\a_m}}$ in this case are known as \emph{top connection coefficients}, and the
above formula was originally given by Goulden and Jackson~\cite{gj:cacti}.

Secondly, we note that Poulalhon and Schaeffer actually define $
    P_a(\mathbf{x})
    =
    \sum_{\l \ptn a }
    \Dop(e_{2\l+1}(\mathbf{x})/\aut{\l}),
$ and ignore the condition $\len {\l} \leq b$ in our definition of $\PGS{a}{b}$. However, replacing
$\smash{\PGS{q}{n-1}}$ with $\smash{P_{q}}$ in Theorem~\ref{thm:ps} has nil effect, since for
$\Dop(e_{2\l+1}(\mathbf{x}))|_{\mathbf{x}=\mathbf{r-2p}}$  to be nonzero some monomial in $e_{2\l+1}(\mathbf{x})$
must be of the form $x_1^{j_1}\cdots x_m^{j_m}$ with $j_i \leq r_i-2p_i$ for all $i$. This implies  $|2\l+1| =
\sum_i j_i \leq \sum_i (r_i-2p_i) =2g+n-1 - \sum_i 2p_i$, while the conditions $\l \ptn q$ and $q+\sum_i p_i=g$
give $|2\l+1|=2q+\len{\l} = 2g - \sum_i 2p_i + \len{\l}$. Thus we require $\len{\l} \leq n-1$ for nonzero
contributions to $P_{q}(\mathbf{r-2p})$.

\begin{lem*} \label{lem:extraction}
Let $s, t_1,\ldots,t_m$ be nonnegative integers. Set $\mathbf{x}=(x_1,\ldots,x_m)$ and
$\mathbf{t}=(t_1,\ldots,t_m)$, and let $f(\mathbf{x})$ be a homogeneous polynomial of total degree
$t_1+\cdots+t_m-s$. Then
$$
    \sbr{\frac{\mathbf{x}^\mathbf{t}}{\mathbf{t}!}}
    (x_1+\cdots+x_m)^{s} f(\mathbf{x}) =
    s!\,\Dop\!\!\left.\big(f(\mathbf{x})\big)\right|_{x_1=t_1,\ldots,x_m=t_m}.
$$
\end{lem*}

\begin{proof}
Consider the case where $f(\mathbf{x})=x_1^{j_1}\cdots x_m^{j_m}$ with $\sum_i j_i=\sum_i t_i-s$. Here
\begin{align*}
    [\mathbf{x}^\mathbf{t}] (x_1+\cdots+x_m)^s f(\mathbf{x})
    &= [\mathbf{x}^\mathbf{t}]
        \sum_{i_1 + \cdots + i_m = s} \frac{s!}{i_1!\,\cdots\,i_m!} x_1^{i_1+j_1} \cdots x_m^{i_m+j_m} \\
    &= \begin{cases}
            \D \frac{s!}{\mathbf{t}!}\, \prod_{i=1}^m \frac{t_i!}{(t_i-j_i)!}
            &\text{if $j_i \leq t_i$ for all $i$}, \\
            0 &\text{otherwise}
        \end{cases} \\
    &=  s!\,\Dop( f(\mathbf{x}))|_{\mathbf{x}=\mathbf{t}}.
\end{align*}
The general result now follows by linearity.
\end{proof}

\noindent\textbf{Proof of Theorem~\ref{thm:ps}}: Comparing~\eqref{eq:RGS} and~\eqref{eq:SGS} we find that, for $\g
=(\g_1,\ldots,\g_l) \ptn n$,
$$
    \RGS{\g}(x,y) = 2^{l-1} \prod_{i=1}^l \g_i \cdot \sum_{p \geq 0} \SGS{p}(\g) x^{n-l-2p} y^{2p+l-1}.
$$
Thus
\begin{equation*}\label{eq:rs}
    \r{\g}_{j,k} =
    \begin{cases}
        \D\frac{2^{\len{\g}-1} \z{\g}}{\aut{\g}}\, \SGS{p}(\l) &\text{if $(j,k)=(n-\len{\g}-2p,\len{\g}+2p-1)$}, \\
        0 &\text{otherwise}.
    \end{cases}
\end{equation*}
From~\eqref{eq:gencx2} and Theorem~\ref{thm:main} we immediately have
\begin{equation*}
    \cx{(n)}{\a_1,\ldots,\a_m}
    = \frac{n^{m-1}}{2^{2g} \prod_{i} \aut{\a_i}}
        \sum_{\mathbf{p}}
        \, \sbr{\frac{\mathbf{x}^{\mathbf{r-2p}}}{\mathbf{(r-2p)}!}}
        \sum_{\len{\l}=n-1} \frac{e_{2\l-1}(\mathbf{x})}{\aut{\l}}
        \prod_{i=1}^m (\len{\a_i}+2p_i-1)! \,\SGS{p_i}(\a_i),
\end{equation*}
where the outer sum extends over all tuples $\mathbf{p}=(p_1,\ldots,p_m)$ of nonnegative integers. Now
\begin{align}
    \sbr{\frac{\mathbf{x}^{\mathbf{r-2p}}}{\mathbf{(r-2p)}!}}
         \sum_{\len{\l}=n-1} \frac{e_{2\l-1}(\mathbf{x})}{\aut{\l}}
    &=  \sum_{s=0}^{n-1} \sum_{\substack{\l \ptn q \\ \len{\l}=n-1-s}}
        \sbr{\frac{\mathbf{x}^{\mathbf{r-2p}}}{\mathbf{(r-2p)}!}}
        \frac{e_1(\mathbf{x})^s}{s!}
        \frac{e_{2\l+1}(\mathbf{x})}{\aut{\l}}, \label{eq:final}
\end{align}
where $q$ is chosen to make $e_1(\mathbf{x})^s e_{2\l+1}(\mathbf{x})$ of total degree $\sum_i (r_i-2p_i)$.  In
particular, if $\l \ptn q$ and $\len{\l}=n-1-s$, then  $e_1(\mathbf{x})^s e_{2\l+1}(\mathbf{x})$ is of degree
$|2\l+1| +s = 2|\l|+\len{\l}+s=2q+n-1$, so we require
$$
    2q+n-1 = \textstyle \sum_i (r_i-2p_i) = (2g+n-1)-\sum_i 2p_i,
$$
or simply $q+p_1 + \cdots + p_m = g$.  Finally, applying the lemma to the RHS of~\eqref{eq:final} results in
$$
    \sum_{s=0}^{n-1} \sum_{\substack{\l \ptn q \\ \len{\l}=n-1-s}}\!\!\!\!\!\!\!
    \Dop\!\!\left.\bigg(\frac{e_{2\l+1}(\mathbf{x})}{\aut{\l}}\bigg)\right|_{\mathbf{x}=\mathbf{r-2p}}
    =
    \PGS{q}{n-1}(\mathbf{r-2p})
$$
and this completes the proof. \qed

\bibliographystyle{amsplain}
\bibliography{bibliography}

\providecommand{\bysame}{\leavevmode\hbox to3em{\hrulefill}\thinspace}
\providecommand{\MR}{\relax\ifhmode\unskip\space\fi MR }
\providecommand{\MRhref}[2]{%
  \href{http://www.ams.org/mathscinet-getitem?mr=#1}{#2}
}
\providecommand{\href}[2]{#2}
\begin{thebibliography}{1}

\bibitem{biane:largecycles}
P.~Biane, \emph{Nombre de factorisations d'un grand cycle.}, S{\'e}m. Lothar.
  de Combinatoire \textbf{51} (2004).

\bibitem{melou-schaeffer:constellations}
M.~Bousquet-M{\'e}lou and G.~Schaeffer, \emph{Enumeration of planar
  constellations}, Adv. in Appl. Math. \textbf{24} (2000), 337--368.

\bibitem{gj:cacti}
I.P. Goulden and D.M. Jackson, \emph{The combinatorial relationship between
  trees, cacti and certain connection coefficients for the symmetric group},
  European J. Combin. \textbf{13} (1992), 357--365.

\bibitem{goupil-schaeffer:largecycles}
A.~Goupil and G.~Schaeffer, \emph{Factoring {$n$}-cycles and counting maps of
  given genus}, European J. Combin. \textbf{19} (1998), 819--834.

\bibitem{bk:macdonald}
I.G. Macdonald, \emph{{S}ymmetric {F}unctions and {H}all {P}olynomials}, Oxford
  University Press (Clarendon), Oxford, 1979.

\bibitem{poulalhon-schaeffer:largecycles}
D.~Poulalhon and G.~Schaeffer, \emph{Factorizations of large cycles in the
  symmetric group}, Discrete Math. \textbf{254} (2002), 433--458.

\bibitem{bk:stanley2}
R.P. Stanley, \emph{{E}numerative {C}ombinatorics: {V}olume 2}, Cambridge
  University Press, 1999.

\end{thebibliography}

\end{document}